\documentclass{article}

\newcommand{\gnl}{^\lceil}
\newcommand{\gnr}{^\rceil}
\newcommand{\sigx}{\sigma(x)}
\newcommand{\subxx}{Sub(x,Num(x))}
\newcommand{\ssubxx}{\mathrm{S}\subxx}

\newcommand{\bfb}[1]{\mbox{\boldmath $ #1 $}}

\begin{document}

\title{Concerning the Representability of Self-Reference in Arithmetic}
\author{Paul Daniel Carr \\ \textit{University of North Carolina at Chapel Hill}} 
\maketitle

\begin{abstract}
Terms in arithmetic of the form $s$ in the formula $s=t(\langle s \rangle)$, with $t$ a term with one free variable and $\langle s \rangle$ denoting the G\"odel number that encodes $s$, are examined by writing the explicit definition of the encoding functions whose representation they include. This is first done with a specific encoding function and system of encoding and then examined more generally. The surprising result of each such construction, involving conventionally defined substitution or diagonalization functions and using conventional systems of encoding, is shown to be a non-terminating symbolic expression. 
\end{abstract}

\section{Introduction and Notation}\label{intro}
This note concerns  certain ``self-referential'' terms in arithmetic.  A clear, early example of these can be found in Hilbert \& Bernays' \textit{Grundlagen der Mathematik}\footnote{as translated in \cite{Priest}}: for any term $t(x)$ with one free variable $x$, there is a closed term $s$ such that the formula $s=t(\langle s \rangle)$ is provable, where $\langle s \rangle$ is the numeral that denotes the G\"odel number of the term $s$. 

The principal matter investigated here is how such a term, which contains a term which represents the diagonalization function referred to by Hilbert and Bernays, would be written using the symbols of a conventional language of the standard model of arithmetic.  This matter is often considered to be unimportant because the diagonalization function can be shown to be defined in terms of the basic functions of zero, identity and succession by the means of composition, primitive recursion and minimization. If the diagonalization function is therefore defined in terms of simple, representable functions, then it is itself representable.  However, in the case of those terms $s$ for which $s=t(\langle s \rangle)$ is supposed to be provable there is a problem regarding the substitution of certain values.  

The first specific case examined in this note takes $t(x)$ to simply be the successor function. This is because the successor symbol is conventionally treated as a primitive symbol, not definable in terms of other more primitive symbols.  In this case, the result may be understood to mean that there is a closed term $s$ such that  $s=\mathrm{S}(\langle s \rangle)$ is provable, where S is the function symbol representing succession.   It is the resulting closed term in arithmetic that denotes the successor to the G\"odel number which encodes that term itself which will be investigated by closely examining how such a term would be written using conventionally defined encoding and substitution functions.  In so doing it will be shown that, for standard systems of encoding, this term cannot be composed of a finite sequence of symbols.  

It will then be clear that the successor symbol is simply one specific example of a symbol other than certain numerals contained in a term which denotes the successor to the G\"odel number that encodes the term itself.  Because there is at least one such symbol in the term which represents the given functions in the cases examined, this denotation cannot be accomplished in a term written using a finite number of symbols.

In the discussion of this observation in arithmetic which follows, the encoding function defined in \cite{Shoenfield} will be used along with definitions and usage from the same source unless otherwise indicated:  the terms \textbf{0, S0, SS0, ...} are the numerals, but for convenience a numeral with $n$ occurrences of \textbf{S} will be referred to with the shorthand notation $\mathbf{k_n}$. 

Expressions comprised of symbols of the language of arithmetic will be written in boldface, while functions and individuals of the standard model or interpretation \textit{N} of arithmetic will be written in italics.  Any sequence of symbols of the language enclosed in half-brackets $\gnl \: \gnr$, rather than angled brackets $\langle \: \rangle$ will be taken to mean the G\"odel number encoding that sequence.  The angled brackets will be reserved for a particular function, to be defined in a later section.

For example, let $f(x)$ indicate an open, unary function in the standard model of arithmetic. Then  $\mathbf{f(x)}$ indicates the representing term, written using the corresponding function symbol if one is defined, of the basic language of arithmetic. $\gnl \mathbf{f(x)}\hspace{.1mm} \gnr$ will indicate the G\"odel number that encodes the sequence $\mathbf{f(x)}$ and $\mathbf{k_{\gnl f(x)\hspace{.1mm} \gnr}}$ will indicate the numeral denoting this G\"odel number. The term being examined here will be referred to with the symbol \boldmath $\sigma$\unboldmath, indicating a particular term of the language, so the specific formula from \cite{Hilbert&Bernays} mentioned above being examined will be written as \boldmath $\sigma=\mathrm{Sk}_{\gnl \sigma \gnr}$ \unboldmath.

\section{}\label{2}

In order to fully examine the formula from \cite{Hilbert&Bernays}, a preliminary observation must first be made.  This observation concerns the fact that in the encoding method used here (as well as in many other encoding methods) no number $p$ can encode the sequence of symbols comprised of the symbol \textbf{S} followed by the numeral $\mathbf{k_p}$ that denotes itself.  This can be shown as follows:

Let $\alpha$ be the number which encodes the symbol \textbf{S}; $\mathbf{k_{\gnl S \gnr}}$ is $\mathbf{k_{\alpha}}$. Then if $p$ encodes an \textbf{S}-symbol as well as its own numeral $\mathbf{k_p}$, then in the encoding system used in \cite{Shoenfield}, as well as prime-factor-based encoding systems, $p$ is at least $\alpha$.  This implies that $\mathbf{k_p}$ has at least as many \textbf{S}-symbols as  $\mathbf{k_{\alpha}}$.  But then $p$ must now at least encode $\alpha$-instances of \textbf{S}. This implies that $\mathbf{k_p}$ has at least as many \textbf{S}-symbols as $\mathbf{k_{\gnl k_{\alpha}\, \gnr}}$. 

This in turn implies that $p$ at least encodes $\mathbf{k_{\gnl k_{\alpha}\, \gnr}}$, so $\mathbf{k_p}$ has at least as many \textbf{S}-symbols as $\mathbf{k_{\gnl k_{\gnl k_{\alpha}\, \gnr}\,\gnr}}$, and so on. In general, let systems of encoding for which it is provable that
\begin{equation}\label{numeralproof}
\mathbf{k_{\gnl S\,\gnr} < k_{\gnl k_{\gnl S\,\gnr}\,\gnr} < k_{\gnl k_{\gnl k_{\gnl S\,\gnr}\, \gnr}\,\gnr} < ...} 
\end{equation}
for any finite number of steps be called \textit{regular}, similarly (but not equivalently) to \cite{Heck}. Thus
\newtheorem{lem}{Lemma}
\begin{lem}\label{Lemma}
for regular systems of encoding, no numeral consisting of a finite number of \textbf{S}-symbols followed by \textbf{0} can denote a number that encodes an additional \textbf{S} as well as itself.
\end{lem}

For regular systems of encoding then no finitary proof written using the basic symbols of arithmetic can contain such a numeral.  For example,
\begin{equation}\label{numeralproof2}
\neg \vdash_{N}  \mathbf{k_p=k_{\gnl Sk_p\,\gnr}}
\end{equation}
which in turn implies that, for $q = Sp$,
\begin{equation}\label{numeralproof4}
\neg \vdash_{N}  \mathbf{k_q=Sk_{\gnl k_q\,\gnr}}
\end{equation}
holds for encoding systems in which (\ref{numeralproof}) is true. This applies to formulas as well:
\begin{equation}\label{numeralproof5}
\neg \vdash_{N} \mathbf{\exists y (y=Sk_n)\:\mbox{  \textit{if}  }\:n=\gnl Sk_{n}\,\gnr} 
\end{equation}
and
\begin{equation}\label{numeralproof6}
\neg \vdash_{N} \mathbf{\exists y (y=Sk_n)\:\mbox{  \textit{if}  }\:n=\gnl\!\exists y (y=Sk_{n})\,\gnr} 
\end{equation}

Thus, for regular systems of encoding the formula that asserts that there exists a successor to the G\"odel number of the formula itself cannot be finitarily proven in arithmetic if the numeral denoting this number is part of the formula. This formula is the result of substituting the G\"odel number of the resulting closed formula for the variable \textbf{x} in the formula $\mathbf{\exists y (y=Sx)}$.  It is clearly not a finite process to do so, but it will be shown that this is precisely what the construction of the term \boldmath $\sigma$ \unboldmath in \boldmath $\sigma=\mathrm{Sk}_{\gnl \sigma \gnr}$ \unboldmath involves.

\section{}\label{3}

In this section the defined functions and encoding method of \cite{Shoenfield} will be used directly as a specific example.  However, it will be argued in section 5 that the result still holds even if a different choice of encoding method or of syntactical arrangement is made.  

Consider the recursive n-ary (n-place) function defined in \cite{Shoenfield} that has as its value the G\"odel number which encodes a sequence of numbers:

\[
\mu x(\beta(x,0)=n\&\beta(x,1)=a_1\&\ldots\&\beta(x,n)=a_n)
\]

The value of this function is the least number $x$ such that the conjuncts included within the parentheses are true.  $a_1 ,\ldots,a_n$ are the numbers, in sequence, which are encoded by the G\"odel number thus specified, and $\beta(x,i)$ is the binary recursive function such that $\beta(x,i)\leq x \stackrel{_{\displaystyle{.}}}{-} 1$ and such that for any sequence of numbers $a_0 ,a_1 ,\ldots,a_{n-1}$ there is a number $x$ such that $\beta(x,i)=a_i$ for all $i<n$. 

It is here taken to be the case that $\beta(x,i)$ is representable; let the term $\mathbf{b}$ with the variables $\mathbf{x,i}$ be taken to be the term representing it.  There will therefore be a  G\"odel number $\gnl \mathbf{b(x,i)} \gnr$ that encodes the term $\mathbf{b(x,i)}$.  The result being demonstrated here is only strengthened if the explicit definition of $\beta(x,i)$ is represented in full wherever it appears.

Now if the above function $\mu x$ is representable as defined in \cite{Shoenfield}, then there is some term $\mathbf{m_{x_1,...,x_n}}$ such that
\[
\vdash _{N} \mathbf{m_{x_1,...,x_n}[k_{a_1},...,k_{a_n}] = k_p} 
\]
if 
\[
\mu x(\beta(x,0)=n\&\beta(x,1)=a_1\&\ldots\&\beta(x,n)=a_n)=p
\]

It is essential to this demonstration to observe that each numeral $\mathbf{k_{a_i}}$ is a sequence of symbols present in the term $\mathbf{m_{x_1,...,x_n}[k_{a_1},...,k_{a_n}]}$.  Each  $\mathbf{k_{a_i}}$ denotes one of the numbers $a_i$ of the sequence of numbers encoded in this way by the G\"odel number $p$, which in turn is denoted by the closed term $\mathbf{m_{x_1,...,x_n}[k_{a_1},...,k_{a_n}]}$.  

If the closed term $\mathbf{m_{x_1,...,x_n}[k_{a_1},...,k_{a_n}]}$ were to denote the G\"odel number which encoded this term itself, then among the symbol-sequences encoded must be the $\mathbf{k_{a_i}}$, which implies that the numbers $a_0,...,a_n$ encode (at least, among other symbols) the very numerals that denote these numbers themselves. This presents various problems related to the preliminary observation made above. In order to make this connection evident, consider the following function $\sigma$ defined as 
\begin{equation}\label{sigmadef1}
\sigma \equiv \mathrm{S}\mu x(\beta(x,0)=n\&\beta(x,1)=a_1\&\ldots\&\beta(x,n)=a_n)
\end{equation}
which is just the composition of the successor function and the function defining the G\"odel number of a sequence of numbers $a_1 ,\ldots,a_n$. Assuming the successor function is represented using the function symbol \textbf{S}, the term representing $\sigma$, which is referred to here by the symbol \boldmath $\sigma$\unboldmath, is defined in terms of \textbf{S} and \textbf{m} such that
\begin{equation}\label{sigmaterm1}
\vdash_{N} \mathbf{Sm_{x_1,...,x_n}[k_{a_1},...,k_{a_n}] = Sk_p}
\end{equation}
if, as before,
\[
\mu x(\beta(x,0)=n\&\beta(x,1)=a_1\&\ldots\&\beta(x,n)=a_n)=p
\]
but where $a_1,\ldots, a_n$ are now the G\"odel numbers which encode the symbols in the term \boldmath $\sigma$\unboldmath.  Rewriting (\ref{sigmaterm1}) using the symbol \boldmath $\sigma$\unboldmath, it is clearly an $n$-ary form of the result from \cite{Hilbert&Bernays} mentioned above:
\[
\vdash_{N} \bfb{\sigma=\mathbf{Sk}_{\gnl \sigma \gnr}}
\]

If the term \boldmath $\sigma$ \unboldmath in this way denotes the successor to the number which encodes itself, the first symbol in the term to be encoded is \textbf{S}, so $\mathbf{k_{a_1}}$ in (\ref{sigmaterm1}) is a numeral denoting a G\"odel number which encodes at least \textbf{S}, but not all of itself as well, due to the observation mentioned above. 

This implies that among the symbols in \boldmath $\sigma$ \unboldmath left to be encoded by subsequent $a_i,\ldots,a_j$, of which there must therefore be at least one,  are the remaining symbols not encoded by the number denoted by $\mathbf{k_{a_1}}$.  However, the same observation applies to $\mathbf{k_{a_i},\ldots,k_{a_j}}$; the numbers that these numerals denote cannot encode them in their entirety as well as the symbols not yet encoded.  This in turn implies that further $a_l,...,a_m$, denoted by  $\mathbf{k_{a_l},\ldots,k_{a_m}}$ (again at least one) encode at least the symbols in $\mathbf{k_{a_i},\ldots,k_{a_j}}$ not yet encoded by the numbers $a_i,...,a_j$, but not all of themselves in their entirety as well, and so on.

This result also applies to a \boldmath $\sigma$ \unboldmath where the succession function is not represented by \textbf{S}, but by some other symbol or sequence of symbols (such as $\mathbf{...+1}$) at some other position in the sequence of \boldmath $\sigma$ \unboldmath than the beginning.  Since it is still part of the term \boldmath $\sigma$\unboldmath, there must be at least one other sequence of symbols $\mathbf{k_i}$ which denotes a number which at least encodes the sequence representing succession.  In this case the comments above regarding $\mathbf{k_1}$ apply to this $\mathbf{k_i}$; it cannot denote a number which encodes both the sequence representing succession as well as all of itself. There must be a another sequence $\mathbf{k_j}$ which denotes a number which at least encodes the remaining symbols in $\mathbf{k_i}$, but not all of itself, and so on. The following has therefore been shown:

\newtheorem{thm}{Theorem}
\begin{thm}\label{theorem1}
For regular systems of encoding, an $n$-ary term which represents (as defined in \cite{Shoenfield}) the function $\mathrm{S}\mu x(\beta(x,0)=n\&\beta(x,1)=a_1\&\ldots\&\beta(x,n)=a_n)$ and which contains other symbols than the numerals denoting the encoded numbers cannot denote the successor to the number which encodes itself within a finite string of symbols.
\end{thm}

The discussion so far has relied heavily upon the observation that no numeral can denote a number that encodes an additional symbol such as the \textbf{S}-symbol as well as the numeral itself. Of course, it is also the case for the encoding function already referred to that a numeral cannot denote a number that encodes just the numeral itself, implying that there is an even stronger result which can be established for this particular encoding function. 

In either case, the \textbf{S}-symbol serves as a convenient stand-in for any other symbols that the term $\mathbf{m}$ contains other than the numerals $\mathbf{k_{a_1},...,k_{a_{n}}}$.  Since this term cannot be constructed without the presence of at least one addtional (function) symbol apart from the numerals $\mathbf{k_{a_1},...,k_{a_{n}}}$, there is at least one such symbol.   

\section{}\label{4}

Here the use of the function $Sub$ in defining the unary function $\sigma(x)$ in \cite{Shoenfield}will be examined. Because the function $\sigma(x)$ is a one-variable function and not a formula, it is not necessary to use the full definition of $Sub(a,b,c)$ applicable to one- or two- variable functions or formulas as found in \cite{Shoenfield}.  For the present discussion the definition will be modified to allow only for unary functions of a given variable:
\begin{eqnarray*}
Sub(x,Num(x)) &\equiv& Num(x)~~\mathrm{if}~Vble(x) \\
                          &\equiv&\langle(a)_0,Sub((a)_1,Num(x))\rangle \\
                          & &~~~~\mathrm{if}~x=\langle(a)_0,(a)_1\rangle \\
                          &\equiv& x~~\mathrm{otherwise}
\end{eqnarray*}
where $Num(x)$ and $Vble(x)$ are defined as in \cite{Shoenfield} and the bracket notation for the encoding function has been used in the 2nd case. In terms of the $\mu x$-function defined above, this definitional case is
\begin{eqnarray*}
\subxx \equiv \mu z(\beta(z,0)=2\&\beta(z,1)=(a)_0\&\ldots \\
\ldots\beta(z,2)=Sub((a)_1,Num(x))) \\
\mathrm{if}~~x=\mu z^\prime(\beta(z^\prime,0)=2\&\beta(z^\prime,1)=(a)_0\&\beta(z^\prime,2)=(a)_1)
\end{eqnarray*}
where $(a)_0$ and $(a)_1$ are the G\"odel numbers of symbols or expressions and $z^\prime$ simply indicates another variable distinct from $x,z$.

Because the function $\subxx$ is defined by cases, any term $\mathbf{m}$ which supposedly represents $\subxx$ cannot be directly encoded (a unique expression number computed for it) in the same manner as a term defined by a single expression. However, even if the cases are joined by logical conjunct symbols ($\&$) and the implications are made fully symbolic (i.e. written in terms of $\neg$ and $\vee$) an expression number cannot be computed or assigned for all the symbols in the explicit definition of $\subxx$  because the definition of $Sub$ includes an instance of $Sub$; it is not an example of a function that is defined only in terms of previously defined functions.

Of course, $Sub$ may be considered to be an example of a function defined by primitive recursion for finite $x$, since $\langle (a)_0,(a)_1 \rangle > (a)_0,(a)_1$ by virtue of the definition of the $\beta$-function.  It is assumed in this kind of definition that the process of defining one $Sub$ with a given argument in terms of another $Sub$ with a smaller argument will eventually terminate becase each succeeding instance of $Sub$ takes as its argument a smaller number than the previous instance. 

This assumption can be seen to hold when $Sub$ takes a definite number as an argument, since there are only a finite number of numbers smaller than the argument, but what about when $Sub$ is given the putative number encoding the term representing itself as its own argument?  The number of times further $\subxx$'s are employed in a given explicit definition of $Sub(a,Num(a))$ depend on what number is substituted for $a$. If each use of $Sub(x,Num(x))$ affects how the overall explicit definition of $Sub(a,Num(a))$ is written, and therefore represented, then it does not appear to be as straightforward how to represent $Sub(p,Num(p))$; that is, when the G\"odel number $p$ of the term representing $\subxx$ itself is substituted for $x$ in $\subxx$. 

The present discussion aims to pursue an answer to this question indirectly by examining the result when the composed function $\mathrm{S}Sub(x,Num(x))$ is given the G\"odel number of its own representing term, which is assumed to exist, as an argument. 

Let the function $\sigma(x)$ be:
\begin{equation}\label{newsigdef}
\sigma(x) \equiv \mathrm{S}Sub(x,Num(x))
\end{equation}

$\sigma(x)$ is now defined in terms of the right hand side of (\ref{newsigdef}). As such, the G\"odel number  \boldmath $\gnl \sigma(\mathrm{x}) \gnr$ \unboldmath, i.e. of the term representing $\sigma(x)$, is just the number $\gnl \mathbf{S\tilde m_{x}} \gnr$, where $\mathbf{S\tilde m_{x}}$ indicates the term representing $\mathrm{S}Sub(x,Num(x))$.  This notation is chosen to indicate that the form of the term representing $Sub(a,Num(a))$ (or at least the relevant defnitional case of it) will be seen to change depending on what number $a$ is substituted for $x$. The succession part of (\ref{newsigdef}) is again taken to be represented by the \textbf{S}-symbol in accordance with \cite{Shoenfield}. Therefore the $Sub(x,Num(x))$-part is represented entirely by $\mathbf{\tilde m_x}$.

For convenience, the number $\gnl \mathbf{S\tilde m_{x}} \gnr$ will sometimes be referred to as $q$, the corresponding numeral of which is $\mathbf{k_q}$.  To be explicit, it must be assumed that there is at least one free instance of $\mathbf{x}$ in the sequence of symbols comprising the term that represents $\ssubxx$.  Substituting $q$ into $\ssubxx$ yields, according to the second case of the definition of $Sub$ above:
\begin{eqnarray}\label{subqq}
\mathrm{S}Sub(q,Num(q))\equiv \mathrm{S}\mu z(\beta(z,0)=2\&\beta(z,1)=\gnl\!\mathbf{S}\gnr\& \nonumber\\
\beta(z,2)=Sub(\gnl \mathbf{\tilde m_x} \gnr,Num(\gnl \mathbf{S\tilde m_x} \gnr)))
\end{eqnarray}
or, using the simplified bracket notation
\begin{equation}\label{subqq2}
\mathrm{S}Sub(q,Num(q))\equiv \mathrm{S}\langle \gnl \mathbf{S} \gnr,Sub(\gnl \mathbf{\tilde m_x} \gnr,Num(\gnl \mathbf{S\tilde m_x} \gnr)\rangle
\end{equation}

Clearly, the explicit definition of (\ref{subqq}) or (\ref{subqq2}) includes at least part of the explicit definition of the closed term $Sub(q,Num(q))$. Whether it is the whole definition depends upon how a definition by cases is treated.

Since $Sub(q,Num(q))$ is thus defined in terms of the function $\mu z(\beta(z,0)=n \& \beta(x,1)=a_1 \&\beta(x,2)=a_2)$, and if this function is representable as previously discussed, then (\ref{subqq}) or (\ref{subqq2}) is represented by a closed term
\begin{equation}\label{repterm1}
\mathbf{Sm_{x_1,x_2}[k_{\gnl S \gnr},k_{q_1}]}
\end{equation}

where $\mathbf{k_{q_1}}$ is the numeral denoting the value of $Sub(\gnl \mathbf{\tilde m_x}\hspace{.1mm} \gnr,Num(\gnl \mathbf{S\tilde m_x}\hspace{.1mm} \gnr)$. 

This implies that $\mathbf{S\tilde m_x[k_q]}$ is actually a term of the form $\mathbf{Sm_{x_1,x_2}[k_{\gnl S \gnr},k_{q_1}]}$, and therefore $\mathbf{\tilde m_x[k_q]}$ is of the form $\mathbf{m_{x_1,x_2}[k_{\gnl S \gnr},k_{q_1}]}$. If $\mathrm{S}Sub(q,Num(q))=p$ then under the assumptions stated so far it follows that
\begin{equation}\label{repterm2}
\vdash _{N} \mathbf{Sm_{x_1,x_2}[k_{\gnl S \gnr},k_{q_1}] = k_p} 
\end{equation}

Now, the function $\subxx$, as well as the unary term which represents it, are defined such that
\begin{equation}\label{subqq3}
Sub(\gnl  \mathbf{\tilde m_x} \gnr,Num(\gnl \mathbf{S\tilde m_x} \gnr))= \gnl\!\!\mathbf{\tilde m_x[k_{\gnl S\tilde m_x \, \gnr}]}\,\gnr = \gnl\!\! \mathbf{\tilde m_x}[\mathbf{k_q}] \gnr 
\end{equation}
but because of what has just been shown about the form of $\gnl \mathbf{\tilde m_x}[\mathbf{k_q}] \gnr$,
\begin{equation}\label{subqq3.5}
\gnl \mathbf{\tilde m_x}[\mathbf{k_q}] \gnr= \gnl\!\mathbf{m_{x_1,x_2}[k_{\gnl S \gnr},k_{q_1}]} \gnr
\end{equation}
or, in terms of the explicit definition of $Sub$,
\begin{eqnarray}\label{subqq4}
Sub(\gnl  \mathbf{\tilde m_x} \gnr,Num(\gnl \mathbf{S\tilde m_x} \gnr)) \equiv \mu z^\prime(\beta(z^\prime,0)=n\;\& \nonumber\\ \beta(z^\prime,1)=a_1\&\beta(z^\prime,2)=a_2)
\end{eqnarray}
where $a_1$ encodes the first symbol in the closed term $ \mathbf{m_{x_1,x_2}[k_{\gnl S \gnr},k_{q_1}]}$ and $a_2$ encodes the second symbol or iteratively encoded sequence of symbols. Among these symbols (which may be among others depending upon how the definition by cases is treated) is the sequence $\mathbf{k_{\gnl S \gnr}}$, the numeral denoting the number which encodes the symbol $\mathbf{S}$.

Putting (\ref{subqq}) and (\ref{subqq3}) together, the result of the manner in which $Sub$ is defined to proceed is that the result of substituting $q$ into $\subxx$ is 
\[
Sub(q,Num(q))= \gnl\!\mathbf{Sm_{x_1,x_2}[k_{\gnl S \gnr},k_{q_1}]} \gnr 
\]
and therefore
\[
\mathrm{S}Sub(q,Num(q))= \mathrm{S}\gnl\!\mathbf{Sm_{x_1,x_2}[k_{\gnl S \gnr},k_{q_1}]} \gnr =p
\]
or, in order to make the connection to the result from \cite{Hilbert&Bernays} obvious by rewriting (\ref{repterm2}) in terms of the symbol \boldmath $\sigma$ \unboldmath,
\begin{equation}\label{repterm3}
\vdash _{N} \bfb{\mathrm{\sigma = Sk_{\gnl \sigma \gnr}}} 
\end{equation}

The result of what has been shown is now as follows: let it be assumed that $\mathrm{S}Sub(x,Num(x))$ is representable as here described, and the representing term can be encoded using a regular system of encoding by a number $q$, which can then be substituted into $\mathrm{S}Sub(x,Num(x))$. Then the resulting value $p$, which can also be written as $\mathrm{S}\gnl\!\mathbf{Sm_{x_1,x_2}[k_{\gnl S \gnr},k_{q_1}]} \gnr$ is denoted by the term $\mathbf{Sm_{x_1,x_2}[k_{\gnl S \gnr},k_{q_1}]}$ representing $\mathrm{S}Sub(q,Num(q))$ by including among the symbols of which it is written numerals denoting numbers which encode these very symbols.  However, as discussed in sections 2\&3 above, if it is true that no numeral can denote a number that encodes an additional symbol such as \textbf{S} as well as the numeral itself, then there is no such finite term $\mathbf{Sm_{x_1,x_2}[k_{\gnl S \gnr},k_{q_1}]}$ or finite number $\gnl\!\mathbf{Sm_{x_1,x_2}[k_{\gnl S \gnr},k_{q_1}]} \gnr$. Thus

\begin{thm}\label{theorem2}
For regular systems of encoding, the unary (1-place) term \boldmath $\sigma$ \unboldmath which represents the defined closed function $\mathrm{S}Sub(q,Num(q))$ and which denotes the successor to the Godel number encoding the term itself cannot be composed of a finite string of symbols.
\end{thm}

This implies that $\mathrm{S}Sub(q,Num(q))$ is not representable in the specific manner discussed, nor in the general manner discussed in the previous section. This in turn implies that under these assumptions as long as there is a single additional symbol apart from the instances of numerals denoting the numbers encoded contained in the term representing $Sub(q,Num(q))$, then it cannot denote the number that encodes itself.

\section{}\label{5}

In this section the underlying problem in constructing the term \boldmath $\sigma$ \unboldmath in the formula \boldmath $\sigma=\mathrm{Sk}_{\gnl \sigma \gnr}$ \unboldmath and which contains a term which represents diagonalization will be examined in a more general fashion.  

Consider a list of unary recursive functions, represented by terms of the language and arranged by some ordering process, such as length and precedence in symbols, or G\"odel number.  

\begin{equation}\label{functions}
	\begin{array}{c}
	\mathbf{f_0(x)}\nonumber\\
	\mathbf{f_1(x)}\\
	\mathbf{f_2(x)}\\
	\vdots
	\end{array}
\end{equation}

By the $\mathbf{f_n(x)}$ are indicated the open terms representing unary (1-place) recursive functions of the universe of arithmetic. Recursive functions refer to those functions constructible in the language of arithmetic from the constant zero function, succession, and the identity functions by composition, primitive recursion and minimization.  An example of minimization is the function $\mu x_{x<...}(...x...)$ defined in \cite{Shoenfield} which has as its value the least number $x$ such that $(...x...)$ is true, or $...$ if there is no such $x$. In the case of the $\beta$-function, there is no need for this bound, since the numbers defined by the $\beta$-function exist for finite sequences.  

Next construct the corresponding list of G\"odel numbers of the entries in (\ref{functions}):

\begin{equation}\label{functionnumbers}
	\begin{array}{c}
	\gnl \mathbf{f_0(x)} \gnr\nonumber\\
	\gnl \mathbf{f_1(x)} \gnr \\
	\gnl \mathbf{f_2(x)} \gnr \\
	\vdots
	\end{array}
\end{equation}
and the corresponding list of G\"odel numerals:

\begin{equation}\label{functionnumerals}
	\begin{array}{c}
	\mathbf{k}_{\gnl \mathbf{f_0(x)} \gnr}\nonumber\\
	\mathbf{k}_{\gnl \mathbf{f_1(x)} \gnr} \\
	\mathbf{k}_{\gnl \mathbf{f_2(x)} \gnr} \\
	\vdots
	\end{array}
\end{equation}

Now construct an array by substituting numerals, beginning with $\mathbf{k_0}$, for $\mathbf{x}$ in each term in (\ref{functions}):
\begin{equation}\label{closedfunctions}
	\begin{array}{cccc}
	\mathbf{f_0(k_0)} & \mathbf{f_0(k_1)} &	\mathbf{f_0(k_2)} & \ldots \\
	\mathbf{f_1(k_0)} & \mathbf{f_1(k_1)} & \ldots                 &        \\
	\mathbf{f_2(k_0)} & \ldots	                &                        &        \\
	\vdots                 &                        &                        &        \\
	\end{array}
\end{equation}

Finally, consider the array of G\"odel numbers of each of the terms in (\ref{closedfunctions}):	
\begin{equation}\label{godelnumbers}
	\begin{array}{cccc}
	\mathbf{\gnl f_0(k_0)\hspace{.1mm}\gnr} & \mathbf{\gnl f_0(k_1)\hspace{.1mm}\gnr} & \mathbf{\gnl	f_0(k_2)\hspace{.1mm}\gnr} & \ldots \\
	\mathbf{\gnl f_1(k_0)\hspace{.1mm}\gnr} & \mathbf{\gnl f_1(k_1)\hspace{.1mm}\gnr} & \ldots                 &        \\
	\mathbf{\gnl f_2(k_0)\hspace{.1mm}\gnr} & \ldots	                &                        &        \\
	\vdots                 &                        &                        &        \\
	\end{array}
\end{equation}

By virtue of their construction as G\"odel numbers, each of these numbers is the value of the \textit{n}-ary function that is used in whatever encoding method has been chosen to encode \textit{n} symbols, or groups of symbols in case a form of encoding compression is used.  In either case, the G\"odel number that encodes the sequence of symbols in a given term in (\ref{closedfunctions}) is the value of such a function and is the corresponding entry in (\ref{godelnumbers}).  One example of such a function has been mentioned above:

\[
\mu x(\beta(x,0)=n\&\beta(x,1)=a_1\&\ldots\&\beta(x,n)=a_n)
\]

\noindent meaning the smallest number $x$ such that $\beta(x,0)=n\&\ldots$ where $n$ is the number of symbols, or encoded groups of symbols, in a closed term $\mathbf{f_i(k_j)}$; $a_1$ encodes the first symbol or group of symbols in the expression, $a_2$ the second, etc.  This, as discussed above, is the function which defines the G\"odel number encoding a sequence of symbols in the notation of \cite{Shoenfield}.  

Whatever choice is made for the encoding method, and thus whatever \textit{n}-ary encoding function yields as value the number encoding $n'$ symbols or encoded groups of symbols in the expression $\mathbf{f_i(k_j)}$, if this function is representable, there will be a term $\mathbf{m_{(ij)}}$ with $\mathbf{k_1,...,k_n}$ for which  
\[
\vdash _{N} \mathbf{m_{(ij)x_1,...,x_n}[k_{a_1},...,k_{a_n}] = k_p} 
\]
i.e., which denotes the number $p$ if $p$ is the value of the function in question when $a_1,...,a_n$ are substituted for its variables $x_1,...,x_n$.

Now, for each entry in (\ref{closedfunctions}), reading the symbols in the entry in sequence, construct the corresponding term \textbf{m} with $\mathbf{k_1,...,k_n}$ which denotes the G\"odel number $p$ encoding the symbol numbers of that entry, i.e., denotes the corresponding entry in (\ref{godelnumbers}):

\begin{equation}\label{godelterms}
	\begin{array}{cccc}
	\mathbf{m_{(00)x_1,...x_n}[k_{a_1},...,k_{a_n}]} & \mathbf{m_{(01)x_1,...x_n}[k_{a_1},...,k_{a_n}]} & \ldots &  \\
	\mathbf{m_{(10)x_1,...x_n}[k_{a_1},...,k_{a_n}]} & \mathbf{m_{(11)x_1,...x_n}[k_{a_1},...,k_{a_n}]} & \ldots                 &        \\
	\mathbf{m_{(20)x_1,...x_n}[k_{a_1},...,k_{a_n}]} &		\ldots									&												 &        \\
		\vdots                   &                            &                        &        \\
	\end{array}
\end{equation}

where the subscripts on the $\mathbf{m}$ are for notational purposes and do not necessarily indicate different symbols of the language, and where in each case the number of variables $x_1,...,x_n$ may be different, depending on how many  numbers $a_1,...,a_n$ are being encoded. These terms will be stepwise written by encoding each symbol or group of symbols in each entry of (\ref{closedfunctions}) with a number $a_i$ and then replacing the corresponding $\mathbf{x_i}$ with $\mathbf{k_{a_i}}$ resulting in the corresponding term in (\ref{godelterms}).  These terms cannot be fully specified unless each symbol of the corresponding entry in (\ref{closedfunctions}) is specified. 

It is essential to note that the array (\ref{godelterms}) is comprised of terms each of which denotes a specific G\"odel number by including, in sequence and as part of the term itself, instances of the G\"odel numerals of the symbols or encoded groups of symbols in the sequence being encoded. Use of the $\beta$-function to accomplish this may or may not be part of a particular encoding method, but the general property just described is of primary importance to the result of this discussion: in order to encode a sequence by a particular number, and write a term which denotes that particular number by the means of representing the encoding function, the full and explicit term which does so must include the numerals denoting (or some other encoded reference to) the G\"odel numbers of the symbols in the sequence being encoded. Although this is a difficult property to prove in general\footnote{in addition to being true of the systems already discussed, it is true of the system orginally used by G\"odel.  See Appendix A} for any possible encoding system, it is also difficult to concieve of a primary encoding function\footnote{one not defined in terms of other encoding functions} that does not require that the numbers to be encoded be specified in some form that can be represented.

Now suppose it were possible to define a unary function $\sigx$  with $x$ free and (by cases if necessary) in terms of the encoding function such that when a G\"odel number $\mathbf{\gnl f_n(x) \gnr}$ which is an entry in (\ref{functionnumbers}) is substituted for $x$ in $\sigx$, the resulting value is the successor to the G\"odel number of that function with its own G\"odel numeral in place of the $\mathbf{x}$:

\begin{equation}\label{sigdef}
\sigma(\mathbf{\gnl f_n(x) \gnr}) = \mathrm{S}\mathbf{\gnl f_n(k_{\gnl f_n(x)\gnr}) \gnr}
\end{equation}

In other words, when a G\"odel number $\mathbf{\gnl f_n(x) \gnr}$ which is an entry in (\ref{functionnumbers}) is substituted for $x$ in $\sigx$, the resulting value is the successor to the corresponding entry  $\gnl \mathbf{f_n(k_i)} \gnr$ in (\ref{godelnumbers}) where $\mathbf{k_i}$ is the numeral $\mathbf{k_{\gnl f_n(x) \gnr}}$.

As an example, using the $\beta$-function method with this substitution the explicit definitional case of $\sigma(\mathbf{\gnl f_{n}(x) \gnr})$ would be, 
\[
\mathrm{S}\mu x(\beta(x,0)=n\&\beta(x,1)=a_1\&\ldots\&\beta(x,n)=a_n)
\]
where the $a_n$ are the symbols or encoded groups of symbols in the ``diagonal'' term $\mathbf{f_n(\gnl f_n(x) \gnr)}$.  If this definitional case is representable, then there is some term $\mathbf{m_{x_1,...,x_n}}$ with $\mathbf{k_1,...,k_n}$ such that
\[
\vdash _{N} \mathbf{Sm_{x_1,...,x_n}[k_{a_1},...,k_{a_n}] = k_p} 
\]
if $p$ is the value of the function $\sigma$ in this case, and again assuming that succession is represented with the symbol \textbf{S}.

Therefore, when a G\"odel numeral from (\ref{functionnumbers}) is substituted for $x$ in $\sigma(x)$ the result of the definition by cases of $\sigma(x)$ is an expression  which is represented by a term $\mathbf{Sm_{x_1,...,x_n}[k_{a_1},...,k_{a_n}]}$ which only differs from the corresponding ``diagonal'' entry in (\ref{godelterms}) by an \textbf{S} symbol.

Now it is possible to consider the result if \boldmath$\sigma(\mathrm{x})$\unboldmath, meaning the representing term of $\sigx$, were to be assigned as one of the terms on list (\ref{functions}); it would have to have some ordering-parameter to assign it a location on the list.  Also,  being written in symbols of the language a G\"odel number could be computed to encode the sequence of symbols in \boldmath $\sigma(\mathrm{x}) $\unboldmath; this G\"odel number would be a corresponding entry in the list (\ref{functionnumbers}).  Let it be written $q$, and the corresponding numeral $\mathbf{k_q}$.  As such, let $q$ be substituted for $x$ in $\sigx$:

\begin{equation}\label{sigma(q)}
\sigma(q)=\mathrm{S}\gnl \bfb{\sigma(\mathrm{k_q})}\gnr
\end{equation}

Readers familiar with the general result from \cite{Hilbert&Bernays} mentioned above will note that the function $\sigma$ in (\ref{sigdef}) may be constructed from the successor function and the recursive function referred to in \cite{Hilbert&Bernays} as diagonalisation.  It will also be evident that a term representing (\ref{sigma(q)}) is simply the special case in \cite{Hilbert&Bernays} mentioned above.

In this case, however, since the term $\sigx$ is an entry on the list (\ref{functions}), \boldmath $\sigma(\mathrm{q})$ \unboldmath is then an entry of (\ref{closedfunctions}) in the row corresponding to \boldmath $\sigma(\mathrm{x})$ \unboldmath and the number \boldmath $\gnl \sigma(\mathrm{q})\gnr$ \unboldmath is the corresponding entry in (\ref{godelnumbers}).  This in turn implies that there is a corresponding term in (\ref{godelterms}) which denotes this number.

The central question to the demonstration is then: what is this entry denoting \boldmath $\gnl \sigma(\mathrm{q}) \gnr$ \unboldmath in the list (\ref{godelterms})?  As already stated, the entries in (\ref{godelterms}) are written by reading the symbols in sequence of the corresponding terms in (\ref{closedfunctions}) and encoding them with the $a_1,a_2,...$ in the terms in (\ref{godelterms}).  But in the case of \boldmath $\sigma(\mathrm{q})$\unboldmath, it is both an entry in (\ref{closedfunctions}) as well as, by definition by cases, (apart from the successor symbol) the entry corresponding to itself in (\ref{godelterms}). To see this, note that (\ref{closedfunctions}) becomes

\begin{equation}\label{closedfunctions2}
	\begin{array}{cccc}
	\mathbf{f_0(k_0)} & \ldots &	\mathbf{f_0}(\mathbf{k_i}) & \ldots \\
	\mathbf{f_1(k_0)} & \ldots & \ldots                 &        \\
	\mathbf{f_2(k_0)} & \ldots	                &        \ddots                &        \\
	\vdots &        &        &   \\
	\bfb{\sigma}(\mathbf{k_0}) &  \ldots   &       & \bfb{\sigma}(\mathbf{k_q})   \\
	\vdots &        &        &  \ldots \\
	\end{array}
\end{equation}
where \boldmath $\sigma(\mathrm{k_q})$ \unboldmath is the representing term of the closed function $\sigma(q)$, which is defined by cases to be the composition of succession and an encoding function which yields as value the G\"odel number $\gnl \bfb{\sigma(\mathrm{k_q})}\gnr$.  Thus \boldmath $\sigma(\mathrm{k_q})$ \unboldmath may be written more explicitly as $\mathbf{Sm_{x_1,...,x_n}[k_{a_1},...,k_{a_n}]}$, where $\mathbf{m_{x_1,...,x_n}[k_{a_1},...,k_{a_n}]}$ denotes the number $\gnl \bfb{\sigma(\mathrm{k_q})}\gnr$.

Meanwhile, the array (\ref{godelterms}) of the terms which represent encoding functions becomes
\begin{equation}\label{godelterms2}
	\begin{array}{cccc}
	\mathbf{m_{(00)x_1,...x_n}[k_{a_1},...,k_{a_n}]} & \ldots & \mathbf{m}_{(0\,i)...}... & \ldots \\
	\mathbf{m_{(10)x_1,...x_n}[k_{a_1},...,k_{a_n}]} & \ldots &  &   \\
	\mathbf{m_{(20)x_1,...x_n}[k_{a_1},...,k_{a_n}]} & \ldots	&  &    \\
		\vdots                   &                      &             \ddots             &        \\
		\mathbf{m_{(\sigma 0)x_1,...,x_n}[k_{a_1},...,k_{a_n}]} & \ldots &  & \mathbf{m_{(\sigma \sigma) x_1,...,x_n}[k_{a_1},...,k_{a_n}]} \\
			\vdots                   &                            &                        &        \\
	\end{array}
\end{equation}

The result is that the entry $\mathbf{m_{(\sigma \sigma)}}$ in (\ref{godelterms2}) which denotes \boldmath $\gnl \sigma(\mathbf{k_q}) \gnr$ \unboldmath is the same sequence of symbols as that following the symbol \textbf{S} in the term in (\ref{closedfunctions2}) that is being read off to produce this very entry.  

This is so even though the term $\mathbf{m_{(\sigma \sigma)}}$ in (\ref{godelterms2}) is a closed term formed from substituting for $n$ variables and the term \boldmath $\sigma(\mathrm{k_q})$ \unboldmath is a closed term formed from substituting for one variable.  To reiterate, the definition by cases of $\sigma(x)$ was such that in the case where $x$ is the G\"odel number of a term $\mathbf{f(x)}$ representing an open unary function, $\sigma(x)$ is simply the successor function composed with the encoding function with $n$ arguments which together yields as its value the number $\mathrm{S}\gnl \mathbf{f(\gnl f(x) \gnr)} \gnr$.  And of course the $\mathbf{m_{(ii)}}$ in (\ref{godelterms2}) are just terms representing encoding functions which yield as their values the G\"odel numbers of the terms in (\ref{closedfunctions2}).

As such, the first symbol which must be encoded by the G\"odel number \boldmath $\gnl \sigma(\mathbf{k_q}) \gnr$ \unboldmath denoted by the corresponding term in (\ref{godelterms2}) is the symbol \textbf{S}, since this is the first symbol read off from the corresponding entry \boldmath $\sigma(\mathbf{k_q})$ \unboldmath in (\ref{closedfunctions2}).

So the first symbol to be encoded by $a_1$, denoted by $\mathbf{k_{a_1}}$ in the term in (\ref{godelterms2}) which denotes \boldmath $\gnl \sigma(\mathbf{k_q}) \gnr$ \unboldmath, is \textbf{S}.  But this $\mathbf{k_{a_1}}$ is a sequence of symbols among the next symbols, after \textbf{S}, of the term in (\ref{closedfunctions2}) being read off.  

As observed in the previous section, in prime-number-based encoding methods no numeral $\mathbf{k_i}$ can denote a number which encodes a sequence of symbols that includes an additional (successor) symbol as well as the numeral itself. In general, for any encoding method for which this is true, there must be another number $a_2$ which encodes the symbols not encoded by $a_1$.

It follows in this case that the next sequence of symbols after \textbf{S}, namely the one that includes $\mathbf{k_1}$  is the next sequence of symbols to be encoded by the $a_2,a_3,...$ and so on.  As soon as an $a_i$ is assigned to encode a previous symbol (or group of symbols) being encoded, there is another sequence of symbols, which includes $\mathbf{k_i}$, as of yet to be encoded.  

In general, any additional symbol or symbols apart from the numerals $\mathbf{k_i}$ in the term $\mathbf{m_{\sigma \sigma}}$ above will serve the same purpose as the symbol \textbf{S} did in \boldmath $\mathrm{\sigma(k_q)}$\unboldmath. Since there must be at least one such (function) symbol in order for the open term $\mathbf{m_{(\sigma)x_1,...,x_n}}$ to represent a well-defined function, the closed term
\[
\mathbf{m_{(\sigma \sigma) x_1,...,x_n}[k_{a_1},...,k_{a_n}]}
\]
cannot therefore denote the number that encodes the term itself.

\section{Conclusion}\label{Conclusion}

Concerning the specific case of the formula from \cite{Hilbert&Bernays},
\[
\bfb{\sigma = \mathrm{Sk}_{\gnl \sigma \gnr}}
\]
it has been shown that the term \boldmath $\sigma$ \unboldmath cannot be constructed within a finite sequence of symbols for the conventional systems of encoding discussed here.  

This followed from the observation that, using the conventional systems of encoding specified, the substitution function cannot take the G\"odel number of the term representing this function as its argument if the resulting closed function is to be represented with a finite sequence of symbols.  This is made clear by the presence of other symbols in such a function than the numerals which denote numbers being encoded (which ones depend upon exactly how the function is represented), such as the additional symbol \textbf{S}.

This implies that, under these assumptions with respect to finite proofs, any term or formula which contains this form of arithmetized self-reference cannot be written in a finite number of symbols, and therefore cannot be part of a finite proof:
\[
\neg \vdash_{N} \bfb{...\sigma...}
\]
or
\[
\neg \vdash_{N} \bfb{...\mathrm{k}_{\gnl \sigma \gnr}...}
\]
where, for any finite number $i$ and the for the encoding systems discussed, it has been shown that \boldmath $\mathrm{k_{\gnl \sigma \gnr} > k_i}$ \unboldmath. 

Generalizing the example discussed in section 4, if $q$ is the G\"odel number of the term representing the function $\mathrm{g}[Sub(x,Num(x))]$, where $\mathrm{g[\mathit{x}]}$ is a single-variable recursive function always represented by at least one symbol, it cannot be proven in the language of arithmetic that $\exists y(y=\mathrm{g}[Sub(q,Num(q))])$ because $\mathrm{g}[Sub(q,Num(q))]$ cannot be finitarily represented. In general, any formula that contains it cannot be proven finitarily within the language of arithmetic. Thus, with \boldmath $\sigma $ \unboldmath the representing term of $\mathrm{g}[Sub(q,Num(q))]$,
\begin{equation}\label{conclusion1}
\neg \vdash_{N} \bfb{\exists \mathrm{y (\sigma = y)}}  
\end{equation}

In this paper only a basic theory $N$ of arithmetic has been employed.  Let induction axioms be added to the theory, such as by the use of an universal induction axiom schema, to produce a theory  of Peano Arithmetic $P$.
Even if $\mathrm{S}Sub(x,Num(x))$ is representable, the expression
\[
\forall x \exists y (y=\mathrm{g}[Sub(x,Num(x))])
\]
cannot be proven in $P$ without being inconsistent with the observation made in (\ref{conclusion1}).

This observation recalls a result mentioned in \cite{Nelson}, concerning the predicate $\phi(x)$ there mentioned as having no proof of being true of particular individuals $b$ in a theory T although it could be shown to be true of the individuals 0,1,2,3,etc..  Although the observation made here is concerned with symbolic considerations and encoding systems rather than using the approach used there it is unlikely that this resemblance is purely a coincidence.  The non-terminating nature of the term \boldmath $\sigma$ \unboldmath could perhaps be described as symbolic impredicativity.  By this is meant that the term \boldmath $\sigma$ \unboldmath cannot denote a particular G\"odel number because of the fact that it ``refers'' to (by means of containing numerals that denote numbers that encode) a symbolic sequence that contains the term itself. 

The connection between the observation (\ref{conclusion1}) and the process of diagonalization is similar.  Diagonalization may be described essentially as the process of using a list or array to define or construct an object not on the list or array.  The analogy to symbolic terms then involves the use of an array of terms to define another term not on this array, not by the use of negation, as is the case with other arguments involving diagonalization (\cite{Simmons},\cite{Smullyan}), but by the means of an additional symbol, such as that defining succession. Of course, if at least part of such a term is on an array associated with (used to generate by means of the encoding system) the one used to define the term in question (the ``diagonalized term'') then a problematic circularity arises which leads directly to the results discussed above. The underlying connection to other forms of the diagonal argument then involves the necessity of extra symbols to indicate terms representing unequal numbers, and the connection between logical negation and inequality.

Finally, these results pertain to the representability of a general class of fixed point arguments, and of other formulas involving substitution defined in terms of an encoding function such as $Sub(x,Num(x))$.  These generalizations will be discussed in later work.  

\section{Acknowledgements}\label{Acknowledgements}

The observations which have been made here are the Providential result of an inquiry into the nature of a proof in arithmetic that there exists a successor to any number.  This inquiry began out of a series of conversations with Frank M. Boardman, who postulated that in every fact-functional counting system there must be a largest number, a number which does not have a successor.  It is this unexpected assertion made by Mr. Boardman which provided the foundation for the inquiry.

The motivation to consider this assertion with respect to physical calculations, as well as to examine the logical foundations of the concept of infinity in arithmetic arose from an interest in the frequent reliance upon, and occasionally problematic nature of this concept in the realm of physical theory. 

This paper would not have been possible without the unwavering endorsement and support of Yee Jack Ng, to whom the author is deeply grateful.  

The author also wishes to acknowledge the useful criticism of Ryan Rohm during a very early stage of the work, and of Keith Simmons during a later stage.   Also, advice from Jim Yuill was very helpful and appreciated.

\appendix
\section{}\label{app1}

Here the demonstration is constructed in the orginal syntax of G\"odel [1931] as reproduced in \cite{Heijenoort}. In this translation, the first prime numbers are \textit{assigned} to each of the basic symbols of the language, and a number $x=2^{n_1}\cdot3^{n_2}\cdot\ldots\cdot p^{n_k}_k$ is \textit{associated} with a sequence $n_1,n_2,...,n_k$ of such numbers, which in turn is associated with a corresponding sequence of basic symbols.  For the purpose of brevity the usage to follow will simply be that a number $x$ is \textit{associated} with a sequence of basic symbols.  

To begin with, the function corresponding to $Sub(a,b,c)$ from \cite{Shoenfield} is Function 31, reproduced in (\ref{function31}), defined in terms of Function 30, reproduced in the two lines of (\ref{function30}):

\begin{eqnarray}
Sb_{0}\left( x^{v}_{y} \right) &\equiv& x \nonumber\\
Sb_{k+1}\left( x^{v}_{y} \right) &\equiv& Su [Sb_k (x^v_y)]\scriptstyle{k\,St\,v\,, x \choose y} \label{function30} \\
Sb(x^v_y) &\equiv& Sb_{A(v,x)}(x^v_y) \label{function31}
\end{eqnarray}
where $k\,St\,v, x$ is Function 28, the value of which is $(k+1)$-th place, counted from the right end of the sequence which is associated with the number $x$, at which the symbol $v$ is free.  $A(v,x)$ is Function 29, the value of which is the number of places in which this occurs.

Function 30 is in turn defined in terms of Function 27:

\begin{eqnarray}\label{function27}
Su [ x\scriptstyle {n \choose y} \displaystyle ] \equiv \scriptstyle \varepsilon \displaystyle z \{ z \le [Pr(l(x)+l(y))]^{x+y} \nonumber \\ 
 \& [(Eu,v)u,v \le x  \:\&\:  x=u\ast\! R(n\: Gl \: x)\!\ast\! v \nonumber \\
 \& \: z=u\!*\! y\!*\! v\: \& \:n=l(u)+1] \}
\end{eqnarray}

Function 31 (equation \ref{function31}) can be defined as a single variable function using Function 19, written $Z(n)$, the value of which is the number associated with the sequence of symbols constituting the numeral \textbf{n}.  Hence $Z(x)$ is substituted for $y$:

\begin{eqnarray}\label{Su(x,Z(x))}
Sb(x^v_{Z(x)}) \equiv Sb_{A(v,x)}(x^v_{Z(x)})
\end{eqnarray}

Thus, substituting for $x$ the number $p$ associated with a term containing the free variable $v$, (\ref{Su(x,Z(x))}) is defined to yield as value the number associated with a new term comprised of the sequence of symbols wherein every free instance of the variable $v$ is replaced with the number associated with the sequence of symbols constituting the numeral \textbf{p}. This new number is explicitly defined in (\ref{function27}) by the expression $z=u*y*v$ which occurs in the final instance of $Su [x{n \choose Z(x)}]$ in $Sb(x^v_{Z(x)})$.

As before, the demonstration proceeds by letting $q$ be the number which is associated with the sequence of symbols representing the full definition of (\ref{Su(x,Z(x))}) preceeded by the successor symbol \textbf{S}.  Let \textbf{q} be the numeral corresponding to the number $q$.

When the last instance of $x$ in (\ref{Su(x,Z(x))}) is replaced with $q$, the full definition of $\mathrm{S}Sb(q^v_{Z(q)})$ includes the last instance of $Su[q{n \choose Z(q)}]$ which has as its definition the following expression:

\begin{eqnarray}\label{Su(q,z(q))}
Su \left[ q\scriptstyle {n \choose Z(q)} \right] \equiv \scriptstyle \varepsilon \displaystyle z \{ z \le [Pr(l(q)+l(Z(q)))]^{q+Z(q)} \nonumber \\ 
 \& [(Eu,v)u,v \le q  \:\&\:  q=u\ast\! R(n\: Gl \: q)\!\ast\! v \nonumber \\
 \& \: z=u\!\ast\! Z(q)\!\ast\! v\: \& \:n=l(u)+1] \}
\end{eqnarray}
where the $z$ being defined here is now the number which is associated with the sequence of symbols which consists of \textbf{S} followed by the term representing the full definition of $\mathrm{S}Sb(q{v \atop Z(q)})$, which includes the expression (\ref{Su(q,z(q))}).

Therefore, the $z$ being defined must be a number which is associated with a sequence that at least includes the symbol \textbf{S}. If so, then because the expression defining $z$ is represented by a sequence to which $z$ is associated, the number $z$ must also be associated with the symbols which constitute a numeral which denotes a number which is at least associated with the symbol \textbf{S}, and so on as before.  

As a result the full definition of (\ref{Su(q,z(q))}) becomes, with R(3) the number associated with the successor symbol:
\begin{equation}\label{selfref}
Su \left[ q\scriptstyle {n \choose Z(q)} \right] \equiv \scriptstyle \varepsilon \displaystyle z {\ldots \& z=\mathrm{R}(3)*\ldots *u'*\ldots}
\end{equation}
where among the part of the expression that follows the R(3) is a $u'$ which is associated with this part of the definition itself:
\[
u'=Z(R(3))
\]
so then (\ref{selfref}) becomes
\[
Su \left[ q\scriptstyle {n \choose Z(q)} \right] \equiv \scriptstyle \varepsilon \displaystyle z {\ldots \& z=\mathrm{R}(3)*\ldots *u'*\ldots*u''*\ldots}
\]
where
\[
u''=Z(Z(R(3)))
\]
and so on.

Of course, there are many other symbols in (\ref{selfref}), such that it may seem unnecessary to include a successor symbol in order to establish the demonstration, but as mentioned above the addition of the successor symbol allows the demonstration to hold regardless of whatever redefinition or compression of the other parts of the expression are made. The successor symbol therefore serves the purpose of a basic symbol that cannot be further compressed or redefined as part of another basic expression.

\end{document}